\documentclass[12pt,a4paper,oneside]{amsart}
\usepackage{amsmath,amsthm,amssymb,amsfonts,mathtools}
\usepackage[a4paper,margin=1in]{geometry}
\usepackage{microtype}
\usepackage{hyperref}
\hypersetup{colorlinks=true,citecolor=blue,linkcolor=blue,urlcolor=blue}

\newtheorem{theorem}{Theorem}[section]
\newtheorem{lemma}[theorem]{Lemma}
\newtheorem{proposition}[theorem]{Proposition}
\newtheorem{corollary}[theorem]{Corollary}
\theoremstyle{definition}
\newtheorem{definition}[theorem]{Definition}
\newtheorem{remark}[theorem]{Remark}
\newtheorem{example}[theorem]{Example}

\DeclareMathOperator{\rad}{rad}

\title{The minimal periodicity for integral bases of pure number fields}
\author{Khai-Hoan Nguyen-Dang}
\address{Morningside Center of Mathematics, Chinese Academy of Sciences, No.\ 55, Zhongguancun East Road, Beijing, 100190, China}
\email{khaihoann@gmail.com, nguyen@amss.ac.cn}
\date{\today}

\begin{document}

\begin{abstract}
Fix $n\ge3$. For the pure field $K_a=\Bbb Q(\theta)$ with $\theta^n=a$, where $a\neq \pm1$ is $n$th-power-free, $x^n-a$ is irreducible over $\Bbb Q$, and Hypothesis~\textup{(\textbf{H})} holds, we encode an integral basis in the fixed coordinate $\{1,\theta,\dots,\theta^{n-1}\}$ by its \emph{shape}. We prove a sharp local-to-global principle: for each $p^e\!\parallel n$, the local shape at $p$ is determined by $a\bmod p^{\,e+1}$, and this precision is optimal. Moreover, the global shape is periodic with minimal modulus
\[
M(n)=\prod_{p^e\parallel n}p^{\,e+1}=n\cdot\mathrm{rad}(n),
\]
providing many applications in the understanding integral bases of pure number fields.
\end{abstract}
\maketitle

\section{Introduction}
A classical problem in algebraic number theory is to study the ring of integers $\mathcal{O}_K$ of a number field $K$, in particular for the pure fields $K=\mathbb{Q}(\theta)$ with $\theta^n=a$ (see \cite{Gaal2019Book} for a modern treatment). Throughout we assume that $a$ is $n$th-power-free, $a\neq\pm1$, and that $x^n-a$ is irreducible over $\Bbb Q$ (equivalently, $[K_a:\Bbb Q]=n$). Recently, two complementary phenomena have guided works on this direction.  First, in small degrees (\(3\le n\le 9\)) Ga\'al--Remete \cite{GR17} showed that the \emph{global pattern} of integral bases is \emph{periodic} as \(a\) varies, with period \(n^2\). In particular, if one writes an integral basis with rational coefficients in the \(\{1,\theta,\dots,\theta^{n-1}\}\)-coordinates, then every prime dividing a denominator must divide \(n\) (a global restriction that already foreshadows a local explanation).  This offers a striking, hands‑on picture of regularity in families of pure fields. Secondly, for general \(n\) and squarefree \(a\), Jakhar--Khanduja--Sangwan \cite{JKS21} gave \emph{explicit} integral bases and proved a global period \(n\cdot\prod_{p\mid n}p\).  Their construction is intrinsically local: it splices together \(p\)-integral bases at primes \(p\mid n\) and encodes, for each \(m=1,\dots,n-1\), a $p$‑adic denominator exponent and a small “linear correction”.  The dependence on \(a\) is controlled by a short list of $p$‑adic invariants (e.g. \(v_p(a^{p-1}-1)\)), revealing that only one extra $p$‑adic digit beyond \(p^{e_p}\parallel n\) is ever needed to fix the local shape at \(p\).  

In this paper we answer the question on how much information about \(a\) do the local computations at each \(p\mid n\) really need?  Equivalently, what is the smallest modulus \(M(n)\) such that, whenever \(a\equiv a'\pmod{M(n)}\), the structure of the integral basis of \(\mathcal O_{\Bbb Q(\sqrt[n]{a})}\) is the same as for \(a'\)? We show that the Jakhar--Khanduja--Sangwan period is in fact \emph{minimal}, in a precise sense, for all \(n\ge3\), and we do so in a way that emphasizes the underlying local mechanism. The first novelty of our paper is we unify and formalize these two phenomena by introducing the concept of \emph{shape} (defined precisely in ~\ref{def:shape}) as the sequence, for \(m=1,\dots,n-1\), consisting of (i) the prime‐power denominators supported on \(p\mid n\) and (ii) the corresponding linear corrections recorded modulo those denominators. We then prove that (see Theorem~\ref{thm:minimal-modulus}) the smallest global modulus governing periodicity of the shape is always
\[
\ M(n)=\prod_{p^{e_p}\parallel n} p^{\,e_p+1}\;=\;n\cdot\mathrm{rad}(n)
\]
In particular \(M(n)\mid n^2\) for all \(n\), with equality \(M(n)=n^2\) if and only if \(n\) is squarefree. Thus we both recover the \(n^2\)-period of \cite{GR17} and pinpoint exactly when it is sharp.

\medskip
\subsection*{Method}
Our proof is short and local-to-global and uses only the explicit $p$-local invariants and bases of Jakhar--Khanduja--Sangwan \cite{JKS21} under Hypothesis~(\textbf{H}). Fix a prime power $p^e\parallel n$ and write $v_p$ for the $p$-adic valuation. When $v_p(a)=0$ we set
\[
r_p(a)=v_p(a^{p-1}-1)-1,\qquad d_p(a)=\min\{r_p(a),e\},
\]
and when $v_p(a)>0$ we have $r_p(a)=-1$ (hence $d_p(a)=-1$). In \cite[Theorem~1.6]{JKS21}, the $p$--denominator exponents $k_{p,m}$ are explicit functions of $p,n,m$ and the truncated invariant $d_p(a)=\min\{v_p(a^{p-1}-1)-1,e\}$ (in the unit branch $v_p(a)=0$). Moreover, the correction terms $\beta_m$ are constructed from congruence data modulo $p^{\,k_{p,m}+1}$ (in particular via an inverse of $a$ modulo $p^{\,k_{p,m}+1}$; see the construction preceding \cite[Theorem~1.6]{JKS21}). Hence the reduced classes $\beta_m\bmod p^{\,k_{p,m}}$ are determined by $a\bmod p^{e+1}$, since $k_{p,m}\le e$. Since $d_p(a)$ is determined by $a\bmod p^{e+1}$, the local $p$-shape is determined by $a\bmod p^{e+1}$.

Already among $p$-adic units one can choose
\[
a\equiv 1+p^{e}u\not\equiv 1\pmod{p^{e+1}}
\qquad\text{and}\qquad
a'\equiv 1+p^{e+1}u'\pmod{p^{e+2}},
\]
with $p\nmid uu'$, for which $r_p(a)=e-1$ while $r_p(a')\ge e$, hence $d_p(a)\ne d_p(a')$. Therefore the
$p$-part of the shape changes, so $p^{e+1}$ is the least local modulus at $p$.

By the Chinese Remainder Theorem, the global shape depends only on the residue class of $a$ modulo $\prod_{p\mid n}p^{e_p+1}=n\cdot\mathrm{rad}(n)$. This gives the optimal upper bound, and local sharpness at every $p\mid n$ makes the modulus minimal.

\subsection*{Applications and further consequences}

The point of view developed here---the \emph{local shape} at each \(p\mid n\) and its \emph{global} gluing---yields several immediate applications. Firstly, we can show that one can compute, for any fixed \(n\), a finite lookup table indexed by the classes \(a\bmod M(n)\) recording the \(p\)–denominator exponents and reduced corrections for every \(p\mid n\). Thanks to Theorem~\ref{thm:minimal-modulus}, the table itself can be filled entirely locally using the explicit local recipes in \cite{JKS21}. This offers a uniform and predictable alternative to ad hoc casework used in small degrees.

The second application concerns local tests for the $p$-part of the index
$[\mathcal O_{K_a}:\Bbb Z[\theta]]$ at primes $p\mid n$. When $p\parallel n$ and $v_p(a)=0$, the dichotomy between $a^{p-1}\not\equiv 1\pmod{p^2}$ and $a^{p-1}\equiv 1\pmod{p^2}$ decides whether $p$ divides $[\mathcal O_{K_a}:\Bbb Z[\theta]]$ (equivalently, whether $\Bbb Z[\theta]$ is $p$--maximal). More generally, for $p^e\parallel n$ and $v_p(a)=0$ the relevant invariant is $d_p(a)=\min\{\,v_p(a^{p-1}-1)-1,\,e\,\}$, and Theorem~\ref{thm:minimal-modulus} shows that this (and with it all $p$--denominator exponents and reduced corrections in an explicit basis) is determined by $a\bmod p^{e+1}$ and cannot be recovered from $a\bmod p^e$ alone.

Finally, because the shape is a class function of \(a\bmod M(n)\), families of degree-$n$ pure-field parameters cut out by finitely many local conditions on \(a\bmod M(n)\) admit \emph{explicit counting and density} statements. Using the equidistribution of \(n\)th–power–free integers in arithmetic progressions and then excluding the sparse reducible binomial parameters, one obtains asymptotics and natural densities for any prescribed union of shape classes; in particular, the global density factorises into a product of local proportions attached to residue classes modulo \(p^{e_p+1}\).

\subsection*{Historical remarks and related results}
The study of integral bases and monogeneity goes back to Dedekind and Hensel and is tightly intertwined with the index, discriminant, and the factorization of rational primes in \(\mathcal O_K\). Dedekind’s criterion and index theorem explain how, when \(p\nmid [\mathcal O_K:\mathbb Z[\theta]]\), the prime-ideal factorization of \(p\mathcal O_K\) is read off from the factorization of the minimal polynomial modulo \(p\) (see, e.g., \cite{Dedekind1878}). Hensel developed local methods that relate common index divisors to the splitting pattern of \(p\) \cite{Hensel1908}. 

When Dedekind’s criterion fails (i.e.\ some \(p\mid [\mathcal O_K:\mathbb Z[\theta]]\)), the modern local machine is based on Newton polygons at \(p\), initiated by Ore in the 1920s and refined into a systematic algorithm by Montes and collaborators: one builds (higher-order) \(p\)-adic Newton polygons, reads residual polynomials along principal sides, and from this recovers \(p\)-adic factorization, indices, discriminants, and \(p\)-integral bases \cite{Ore28,MN92,GMN12,GMN15}. This underlies much current algorithmic number theory around integral bases and prime decomposition. For practical $p$-integral basis extraction and index computations using Newton polygons, see ElFadil and his collaborators’ papers \cite{EMN12,Elf20} and expositions cited therein. For the complexity and implementation discussions, take a look at \cite{BNS13}. These provide worked implementations that agree with and conceptually support the $p$-local pieces used by Ga\'al--Remete and Jakhar--Khanduja--Sangwan.

In Kummer or radical extensions over cyclotomic bases, the decisive local congruences at tame primes naturally live one level deeper in the local parameter: modulo \((1-\zeta_p)^2\) (the uniformizer at \(p\) in \(\mathbb Z[\zeta_p]\)), mirroring the pure-field “one more \(p\)-adic digit” from \(p\) to \(p^2\). For a modern treatment proving relative monogeneity criteria and non-monogeneity results in this framework see Smith \cite{Smi19,Smi21}. 

\subsection*{Organization of the paper}
Section~\ref{sec:defs} introduces the \emph{shape} of an integral basis and gives illustrative examples. Section~\ref{sec:proof} is preceded by our main local–to–global result (Theorem~\ref{thm:minimal-modulus}) and its direct consequences. We also include a tame relative version (Proposition~\ref{prop:relative}) and examples. Section~\ref{sec:applications} develops applications: sharp local power-basis, explicit discriminant jumps, and counting/density results for prescribed shape classes, with worked quartic and sextic cases. 

\subsection*{Acknowledgments}
We are grateful to the Morningside Center of Mathematics, Chinese Academy of Sciences, for its support and a stimulating research environment. We really appreciate the anonymous reviewer for their careful reading of the manuscript and for their valuable comments and suggestions. We thank Istv\'an Ga\'al for his interest in our work and for his helpful comments.

\subsection*{Conflict of Interest}
The author declares that there is no conflict of interest.

\subsection*{Data availability}
This research is purely theoretical; no datasets were generated or analyzed. All mathematical results and worked examples are contained in the article.

\subsection*{Funding}

The author received no financial support for the research, authorship, or publication of this article.

\section{Shape, periodicity, and minimal period}\label{sec:defs}

Throughout let $n\ge3$ and $K_a=\Bbb Q(\theta)$ with $\theta^n=a$, where $a\in\Bbb Z$ is $n$th‑power‑free, $a\neq\pm1$, and $x^n-a$ is irreducible over $\Bbb Q$. Write $n=\prod_{p}p^{e_p}$ and $\mathrm{rad}(n)=\prod_{p\mid n}p$.  When we write an integral basis in the $\{1,\theta,\dots,\theta^{n-1}\}$–coordinates we mean a $\Bbb Z$–basis of the form
\[
\Bigl\{\,1,\ \frac{\theta^m+\beta_m}{D_m}\ :\ 1\le m\le n-1\Bigr\},
\]
Here $D_m$ denotes the (normalized) denominator of the $m$th basis element in these fixed coordinates, after clearing the factors supported on primes dividing $a$ into the numerator (cf.\ Remark~\ref{fact:denom}),
with $D_m\in\Bbb Z_{\ge1}$, $\beta_m\in\Bbb Z[\theta]$ and linear in $\{1,\theta,\dots,\theta^{m-1}\}$. We use the following standard assumption in \cite{JKS21} under which the explicit integral bases of $\mathcal O_{K_a}$ are stated. In particular, it holds whenever $a$ is squarefree or $(a,n)=1$. 

\medskip
\noindent
\textbf{Hypothesis (\textbf{H}):} We say that the pair $(n,a)$ satisfies (\textbf{H}) if $a$ is $n$th–power–free and, for every prime $p\mid n$, either $v_p(a)=0$ or the exponent $v_p(a)$ is \emph{not} divisible by $p$
(i.e.\ $p\nmid v_p(a)$).  

Write $n=\prod_{p}p^{e_p}$ and factor $a$ (uniquely up to sign) as
\[
a=\pm\prod_{j=1}^{n-1} a_j^{\,j},
\]
with the $a_j$ squarefree and pairwise coprime. For $0\le m\le n-1$ set
\[
C_m(a)\ :=\ \prod_{j=1}^{n-1} a_j^{\lfloor jm/n\rfloor}.
\]
For each $p\mid n$ define
\[
r_p(a)\ :=\
\begin{cases}
v_p(a^{p-1}-1)-1,& v_p(a)=0,\\[2pt]
-1,& v_p(a)>0,
\end{cases}
\qquad
d_p(a)\ :=\ \min\{\,r_p(a),\,e_p\,\}.
\]
For each $m$ and each $p\mid n$, let $k_{p,m}$ be the largest integer $k$ with $0\le k\le d_p(a)$
such that $m\ge n-n/p^{\,k}$ (put $k_{p,m}=0$ if no such $k$ exists). We obtain the following remarkable result, which is a reformulation of \cite[Theorem~1.6]{JKS21}.
The correction terms $\beta_m$ (and hence their reductions modulo $p^{k_{p,m}}$) are given explicitly by the construction preceding \cite[Theorem~1.6]{JKS21}, based on the auxiliary elements $\delta_{i,m}$ from \cite[Lemma~1.2]{JKS21}; see also \cite[Remark~1.7]{JKS21} for the freedom in choosing inverses.

\begin{theorem}[Explicit local/global form of an integral basis]\label{fact:JKS-form}
Assume (\textbf{H}). Then there exist elements
$\beta_m\in\Bbb Z[\theta]$, linear in $\{1,\theta,\dots,\theta^{m-1}\}$, such that
\[
\mathcal B\ :=\ \Bigl\{\,1,\ \frac{\theta^m+\beta_m}{\;C_m(a)\cdot\prod_{p\mid n}p^{\,k_{p,m}}\;}\ :\ 1\le m\le n-1\Bigr\}
\]
is an integral basis of $\mathcal O_{K_a}$.  Moreover, for each $p\mid n$ the residue class
$\beta_m\bmod p^{\,k_{p,m}}$ is well determined by $p,n,m$ and $a$ via the construction, and
$C_m(a)$ is supported only at primes dividing $a$. 
\end{theorem}

Fix $p\mid n$.  For each $1\le m\le n-1$ let $k_{p,m}\ge0$ be the exponent of $p$ in the denominator of the $m$th basis element when an integral basis is written in the
$\{1,\theta,\dots,\theta^{n-1}\}$–coordinates. In the explicit bases of Theorem~\ref{fact:JKS-form} a linear correction $\beta_m$ appears; we record its reduction modulo the full $p$–power denominator. We formalize the insight in the following notion.

\begin{definition}[Local and global shape]\label{def:shape}
For $p\mid n$, the \emph{local $p$–shape} of $\mathcal O_{K_a}$ is given by
\[
\mathbf S_p(a)\ :=\ \Bigl((k_{p,m})_{m=1}^{n-1},\ (\beta_m \!\!\!\!\pmod{p^{k_{p,m}}})_{m=1}^{n-1}\Bigr).
\]
The \emph{shape} of $\mathcal O_{K_a}$ is the tuple of local shapes 
\[
\mathbf S(a):=(\mathbf S_p(a))_{p\mid n}.
\]
\end{definition}

\begin{remark}[Convention for $k_{p,m}=0$]\label{rem:kpm0}
If $k_{p,m}=0$, then we regard the residue class $\beta_m\bmod p^{\,k_{p,m}}$ as trivial; in particular the $p$–local shape depends only on those $\beta_m$ with $k_{p,m}>0$.
\end{remark}

Following \cite{JKS21}, we \emph{ignore} the factors $C_m(a)$ coming solely from primes dividing $a$.

\begin{remark}[No denominator primes away from $n$ in the normalized form]\label{fact:denom}
Write an integral basis in the $\{1,\theta,\dots,\theta^{n-1}\}$–coordinates in the normalized manner of Definition~\ref{def:shape}, i.e.\ with all factors supported on primes dividing $a$ absorbed into the numerators $C_m(a)$. Then every prime divisor of each denominator $D_m$ divides $n$. Indeed, in the explicit bases of Jakhar--Khanduja--Sangwan the denominators are exactly $\prod_{p\mid n}p^{k_{p,m}}$, while all primes dividing $a$ appear only in $C_m(a)$. Hence, after normalizing by clearing the $C_m(a)$ into the numerators, no prime $q\nmid n$ can appear in any $D_m$. 
\end{remark}

\begin{definition}[Periodicity and minimal period]\label{def:periodicity}
Fix $n$ and let $N\ge1$. We say that the shape is \emph{periodic modulo $N$} if, for all
$n$th--power--free integers $a,a'\in\Bbb Z$ for which $\mathbf S(a)$ and $\mathbf S(a')$ are defined,
\[
a\equiv a'\pmod N\qquad\Longrightarrow\qquad \mathbf S(a)=\mathbf S(a').
\]
A positive integer $N$ with this property is called a \emph{period} (for fixed $n$). The \emph{minimal period}
is the least such $N$.
\end{definition}

We will show that the shape is periodic modulo
\begin{equation}\label{eq:Mn}
M(n)\ :=\ \prod_{p^{e_p}\parallel n} p^{\,e_p+1}\ =\ n\cdot\mathrm{rad}(n),
\end{equation}
and the local sharpness will force $M(n)$ is \emph{minimal}.

\begin{example}[Small degrees $3\le n\le 9$]
Formula \eqref{eq:Mn} gives
\[
M(3)=9,\ M(4)=8,\ M(5)=25,\ M(6)=36,\ M(7)=49,\ M(8)=16,\ M(9)=27.
\]
Thus Ga\'al--Remete’s $n^2$–periodicity for $3\le n\le 9$ is minimal iff $n$ is squarefree; otherwise
$M(n)<n^2$, in agreement with their observation.   
\end{example}

We record two examples that instantiate the definitions of the local and global shape and line up with known tables.

\begin{example}[Quartic case $n=4$]\label{ex:quartic-shape}
Assume Hypothesis (\textbf{H}). Let $a$ be 4th‑power‑free. Put $r=v_2(a-1)-1$ when $a$ is odd, and $r=-1$ when $2\mid a$. Then $d=\min\{r,2\}$. The general rule above gives

$$
k_{2,1}=0,\qquad
k_{2,2}=\mathbf 1_{\{d\ge1\}},\qquad
k_{2,3}=
\begin{cases}
2,& d=2,\\
1,& d=1,\\
0,& d=0.
\end{cases}
$$

So there are exactly three 2‑adic denominator patterns among odd $a$, according to $a\bmod 8$ (and no 2‑denominators when $2\mid a$). Using Ga\'al--Remete’s explicit bases, one can choose the accompanying corrections $\beta_m$ (reduced modulo the full 2‑power denominator) as follows :

\noindent\textbf{Quartic case \(n=4\): basis and corrections (mod \(8\)).}
\begin{center}
\small
\setlength{\tabcolsep}{4pt}
\renewcommand{\arraystretch}{1.15}
\begin{tabular}{@{}p{0.11\linewidth} p{0.17\linewidth} p{0.38\linewidth} p{0.28\linewidth}@{}}
\hline
\textbf{$a\bmod 8$} &
\textbf{$(k_{2,1},k_{2,2},k_{2,3})$} &
\textbf{basis fragment} &
\textbf{reduced corrections}\\
\hline
$3,7$ & $(0,0,0)$ &
$\{\,1,\theta,\theta^2,\theta^3\,\}$ &
$\beta_2\equiv 0,\ \ \beta_3\equiv 0\ (\bmod 1)$\\[2pt]
$5$ & $(0,1,1)$ &
$\bigl\{\,1,\theta,\frac{\theta^2+1}{2},\frac{\theta^3+\theta}{2}\,\bigr\}$ &
$\beta_2\equiv 1\ (\bmod 2),\ \ \beta_3\equiv \theta\ (\bmod 2)$\\[2pt]
$1$ & $(0,1,2)$ &
$\bigl\{\,1,\theta,\frac{\theta^2+1}{2},\frac{\theta^3+\theta+\theta^2+1}{4}\,\bigr\}$ &
$\beta_2\equiv 1\ (\bmod 2),\ \ \beta_3\equiv 1+\theta+\theta^2\ (\bmod 4)$\\[2pt]
\hline
$2,6$ & $(0,0,0)$ &
$\{\,1,\theta,\theta^2,\theta^3\,\}$ &
$\beta_2\equiv 0,\ \ \beta_3\equiv 0\ (\bmod 1)$\\
\hline
\end{tabular}
\end{center}

Local shape $\mathbf S_2(a)$ is completely determined by $a\bmod 8$ (this is the minimal local modulus $2^{\,e_2+1}$). Global shape equals the local one here (since $2$ is the only denominator prime). These three patterns recur with period $2^{e_2+1}=8$, in agreement with both our theorem and the quartic tables in \cite{GR17}. 
\end{example}

\begin{example}[Sextic case $n=6$]\label{ex:Sextic-shape}
Assume Hypothesis (\textbf{H}) (in particular $(a,6)=1$ in the tables below). Here $p=2,3$ with $e_2=e_3=1$ and $M(6)=2^{2}3^{2}=36$. The local moduli are $4$ at $p=2$ and $9$ at $p=3$. We tabulate the \emph{local shapes} and then the \emph{global} CRT gluing.

\noindent\textbf{Local table at \(p=2\) (mod \(4\)).}
\begin{center}
\small
\setlength{\tabcolsep}{4pt}
\renewcommand{\arraystretch}{1.15}
\begin{tabular}{@{}p{0.17\linewidth} p{0.28\linewidth} p{0.19\linewidth} p{0.30\linewidth}@{}}
\hline
\textbf{$a$ (mod $4$)} &
\textbf{$(k_{2,1},k_{2,2},k_{2,3},k_{2,4},k_{2,5})$} &
\textbf{2-part denoms} &
\textbf{reduced $2$-adic $\beta$'s}\\
\hline
$1$ & $(0,0,1,1,1)$ & $m=3,4,5:\ 2$ &
$\beta_3\equiv 1,\ \ \beta_4\equiv \theta,\ \ \beta_5\equiv \theta^2\ (\mathrm{mod}\ 2)$\\[2pt]
$0,2,3$ & $(0,0,0,0,0)$ & none & trivial\\
\hline
\end{tabular}
\end{center}

\noindent\textbf{Local table at \(p=3\) (mod \(9\)).}
\begin{center}
\small
\setlength{\tabcolsep}{4pt}
\renewcommand{\arraystretch}{1.15}
\begin{tabular}{@{}p{0.19\linewidth} p{0.28\linewidth} p{0.17\linewidth} p{0.30\linewidth}@{}}
\hline
\textbf{$a$ (mod $9$)} &
\textbf{$(k_{3,1},k_{3,2},k_{3,3},k_{3,4},k_{3,5})$} &
\textbf{3-part denoms} &
\textbf{reduced $3$-adic $\beta$'s}\\
\hline
$1$ & $(0,0,0,1,1)$ & $m=4,5:\ 3$ &
$\beta_4\equiv 1+\theta^2,\ \ \beta_5\equiv \theta+\theta^3\ (\mathrm{mod}\ 3)$\\[2pt]
$8$ & $(0,0,0,1,1)$ & $m=4,5:\ 3$ &
$\beta_4\equiv 1+2\theta^2,\ \ \beta_5\equiv \theta+2\theta^3\ (\mathrm{mod}\ 3)$\\[2pt]
$0,2,3,4,5,6,7$ & $(0,0,0,0,0)$ & none & trivial\\
\hline
\end{tabular}
\end{center}

\medskip
\noindent\textbf{Global (CRT) gluing modulo \(36\).}
\begin{center}
\small
\setlength{\tabcolsep}{4pt}
\renewcommand{\arraystretch}{1.15}
\begin{tabular}{@{}p{0.40\linewidth} p{0.17\linewidth} p{0.17\linewidth} p{0.20\linewidth}@{}}
\hline
\textbf{Congruence on \(a\)} & \(\mathbf{(k_{2,m})}\) & \(\mathbf{(k_{3,m})}\) & \textbf{Consequence}\\
\hline
\(a\equiv 1\ (\mathrm{mod}\ 4),\ a\equiv \pm1\ (\mathrm{mod}\ 9)\)
&
\(k_{2,3}=k_{2,4}=k_{2,5}=1\)
&
\(k_{3,4}=k_{3,5}=1\)
&
\(m=4,5\) have denom.\ \(2\cdot 3=6\)\\[2pt]
\(a\equiv 1\ (\mathrm{mod}\ 4),\ a\not\equiv \pm1\ (\mathrm{mod}\ 9)\)
&
\(k_{2,3}=k_{2,4}=k_{2,5}=1\)
&
all \(0\)
&
only \(2\)-adic denominators at \(m=3,4,5\)\\[2pt]
\(a\not\equiv 1\ (\mathrm{mod}\ 4),\ a\equiv \pm1\ (\mathrm{mod}\ 9)\)
&
all \(0\)
&
\(k_{3,4}=k_{3,5}=1\)
&
only \(3\)-adic denominators at \(m=4,5\)\\[2pt]
otherwise
&
all \(0\)
&
all \(0\)
&
no denominators\\
\hline
\end{tabular}
\end{center}

All data are determined by $a\ (\mathrm{mod}\ 36)$; the $\beta$‑reductions are taken modulo the full $p$‑power denominator in each row.
\end{example}

\section{Main theorems}\label{sec:proof}
In this section we prove the precise local-to-global statement that underlies the paper: at each prime $p^e\parallel n$ the \emph{local $p$–shape} of $\mathcal O_{K_a}$ is decided by a single additional $p$‑adic digit beyond $p^e$, and this is best possible; consequently the \emph{global} shape is periodic with minimal modulus
\[
M(n)\;=\;\prod_{p^e\parallel n}p^{\,e+1}\;=\;n\cdot\rad(n).
\]
As a preliminary reduction, if an integral basis is written in the $\{1,\theta,\dots,\theta^{n-1}\}$–coordinates, by the Ga\'al--Remete theorem, only primes dividing $n$ can occur in the denominators in the chosen coordinates, so the local shape $\mathbf S_p(a)$ is only defined and only varies for $p\mid n$. 

The proof of our main theorem uses only the explicit local invariants and bases in
Jakhar--Khanduja--Sangwan \cite{JKS21} (under Hypothesis~(\textbf{H})) and a Chinese remainder argument.

We now prove our main result.

\begin{theorem}[Minimal local/global modulus for the shape]\label{thm:minimal-modulus}
Let $n\ge 3$, $a\in\Bbb Z$ be $n$th–power–free with $a\neq\pm1$, and assume $x^n-a$ is irreducible over $\Bbb Q$, and $K_a=\Bbb Q(\theta)$ with $\theta^n=a$ satisfying hypothesis (\textbf{H}). For each prime $p\mid n$ write $p^{e_p}\parallel n$ and let $\mathbf S_p(a)$ and $\mathbf S(a)$ be the local and global shapes.
\begin{enumerate}
\item[\textup{(i)}] \textup{(Local determinacy.)} For every $p\mid n$, the local shape $\mathbf S_p(a)$ depends only on the residue class of $a$ modulo $p^{\,e_p+1}$. Equivalently, if $a'$ is another $n$th--power--free integer satisfying the same standing hypotheses (in particular (\textbf{H}) and irreducibility) and $a\equiv a'\pmod{p^{\,e_p+1}}$, then $\mathbf S_p(a)=\mathbf S_p(a')$.

\item[\textup{(ii)}] \textup{(Local sharpness.)} For every $p\mid n$ there exist $n$th–power–free integers $a,a'$ with $a\equiv a'\pmod{p^{e_p}}$ but $a\not\equiv a'\pmod{p^{e_p+1}}$ such that $\mathbf S_p(a)\ne\mathbf S_p(a')$. Thus $p^{e_p}$ precision does not suffice in general.
\item[\textup{(iii)}] \textup{(Global minimal period.)} The global shape $\mathbf S(a)$ is periodic modulo
\[
M(n)\;:=\;\prod_{p^{e_p}\parallel n} p^{\,e_p+1}=n\cdot\rad(n),
\]
and $M(n)$ is the \emph{least} such modulus: if $N\ge1$ has the property that for all $n$th--power--free $a,a'$ satisfying the standing hypotheses,
\[
a\equiv a'\pmod N\qquad\Longrightarrow\qquad \mathbf S(a)=\mathbf S(a'),
\]
we get $M(n)\mid N$.
\end{enumerate}
\end{theorem}

\begin{proof}
Fix a prime $p\mid n$ with $p^{e}\parallel n$ and write $v=v_p$.
Recall the invariants
\[
r_p(a)=
\begin{cases}
v_p(a^{p-1}-1)-1,& v_p(a)=0,\\
-1,& v_p(a)>0,
\end{cases}
\qquad
d_p(a)=\min\{r_p(a),e\},
\]
and the denominator exponents $k_{p,m}$ defined in \S\ref{sec:defs}.
By \cite[Theorem~1.6]{JKS21}, the exponent \(k_{p,m}\) is determined by
\(p,n,m\) and \(d_p(a)\).  Once \(k_{p,m}\) is fixed, the corresponding
reduced correction term is obtained from the explicit congruence construction
preceding \cite[Theorem~1.6]{JKS21}. In particular, it uses congruence data
modulo \(p^{k_{p,m}+1}\), equivalently an inverse of \(a\) modulo
\(p^{k_{p,m}+1}\) as in \cite[Remark~1.7]{JKS21}.  Thus the correction
\(\beta_m\bmod p^{k_{p,m}}\) is determined by \(a\bmod p^{k_{p,m}+1}\), and
not merely by \(d_p(a)\).

\smallskip
\noindent\textup{(i) Local determinacy.}
Assume $a\equiv a'\pmod{p^{e+1}}$.
Then $v_p(a)>0$ if and only if $v_p(a')>0$.
If $v_p(a)>0$, then $r_p(a)=r_p(a')=-1$ and hence $d_p(a)=d_p(a')=-1$; therefore $k_{p,m}=0$ for all $m$.
By Remark~\ref{rem:kpm0} the reduced correction data are trivial, so $\mathbf S_p(a)=\mathbf S_p(a')$.

Assume now $v_p(a)=0$ (hence also $v_p(a')=0$). From $a\equiv a'\pmod{p^{e+1}}$ we have $a^{p-1}\equiv (a')^{p-1}\pmod{p^{e+1}}$. Therefore
\[
\min\bigl\{v_p(a^{p-1}-1),\,e+1\bigr\}
=
\min\bigl\{v_p((a')^{p-1}-1),\,e+1\bigr\}.
\]
It follows that
\[
d_p(a)=\min\{\,v_p(a^{p-1}-1)-1,\,e\,\}= \min\{\,v_p((a')^{p-1}-1)-1,\,e\,\}=d_p(a').
\]

Since the exponents $k_{p,m}$ are defined purely from $d_p(a)$ (together with $p,n,m$), we get $(k_{p,m})_{m=1}^{n-1}=(k'_{p,m})_{m=1}^{n-1}$.
Finally, in the construction underlying \cite[Theorem~1.6]{JKS21}, each reduced class
$\beta_m \bmod p^{\,k_{p,m}}$ is computed from congruence data modulo $p^{\,k_{p,m}+1}$
(in particular from an inverse of $a$ modulo $p^{\,k_{p,m}+1}$). Since $k_{p,m}\le e$ and
$a\equiv a'\pmod{p^{e+1}}$, the required congruence data agree for $a$ and $a'$, hence the
reduced correction terms coincide. Therefore $\mathbf S_p(a)=\mathbf S_p(a')$.

\smallskip
\noindent\textup{(ii) Local sharpness.}
Take $p$-adic units
\[
a=1+p^{e}u,\qquad a'=1+p^{e+1}u',\qquad p\nmid uu'.
\]
By replacing $a$ and $a'$ by congruent integers if necessary, we may (and do) assume that
$a$ and $a'$ are squarefree and coprime to $n$; in particular they are $n$th--power--free,
satisfy hypothesis (\textbf{H}), and $x^n-a$, $x^n-a'$ are irreducible over $\Bbb Q$.
(For example, each primitive residue class modulo $p^{e+1}$ contains infinitely many squarefree integers.) Then $a\equiv a'\pmod{p^{e}}$ but $a\not\equiv a'\pmod{p^{e+1}}$.
For odd $p$, the lifting-the-exponent lemma gives $v_p(a^{p-1}-1)=e$ and $v_p((a')^{p-1}-1)\ge e+1$; for $p=2$ this is the same statement with $a^{p-1}-1=a-1$.
Thus $r_p(a)=e-1$ while $r_p(a')\ge e$, hence $d_p(a)=e-1$ and $d_p(a')=e$.
Let $m_0:=n-n/p^{e}$. By definition of $k_{p,m}$ we have $k_{p,m_0}=e-1$ but $k'_{p,m_0}=e$, so $\mathbf S_p(a)\ne\mathbf S_p(a')$. This proves \textup{(ii)}.

\smallskip
\noindent\textup{(iii) Global minimal period.}
Part \textup{(i)} shows that for each $p\mid n$ the local shape depends only on $a\bmod p^{e_p+1}$. By the Chinese Remainder Theorem, the global shape depends only on
$a\bmod\prod_{p^{e_p}\parallel n}p^{e_p+1}=M(n)$, so $M(n)$ is a period.
Conversely, if $N$ is any modulus such that $a\equiv a'\pmod N$ implies $\mathbf S(a)=\mathbf S(a')$ for all $n$th–power–free $a,a'$, then for each fixed $p\mid n$ it implies $\mathbf S_p(a)=\mathbf S_p(a')$. 

It remains to justify the divisibility assertion. Suppose, to the contrary, that
\(v_p(N)\le e\) for some \(p^e\parallel n\). Write \(N=p^fN_0\), with
\(f\le e\) and \((N_0,p)=1\). Choose two residue classes modulo a common
modulus divisible by \(N_0\) and \(p^{e+1}\), agreeing modulo \(N_0\), such that
locally at \(p\) they are congruent modulo \(p^e\) but not modulo \(p^{e+1}\),
for instance
\[
A\equiv 1+p^eu\pmod{p^{e+1}},\qquad
A'\equiv 1\pmod{p^{e+1}},
\quad p\nmid u.
\]
By the Chinese Remainder Theorem, we may impose simultaneously
\(A\equiv A'\pmod{N_0}\). Choosing representatives in these classes which
satisfy the standing hypotheses, we get \(a\equiv a'\pmod N\) but
\(\mathbf S_p(a)\ne\mathbf S_p(a')\), contradicting that \(N\) is a period.
Hence \(v_p(N)\ge e+1\).
\end{proof}

\begin{remark}[Comparison with explicit periodicity]
In the explicit construction of integral bases for pure fields, Jakhar--Khanduja--Sangwan show that the $p$–local exponents and corrections are functions of $d_p(a):=\min\{v_p(a^{p-1}-1)-1,e_p\}$ (for $v_p(a)=0$) and of $v_p(a)$, and that $d_p$ is stable under $a\bmod p^{e_p+1}$. This yields periodicity modulo $M(n)=\prod p^{e_p+1}$, and our Theorem~\ref{thm:minimal-modulus} shows that this modulus is also minimal.
\end{remark}

\begin{example}[Nonic fields ($n=9$)]
Here $n=3^2$, so only $p=3$ contributes and the local period is $3^{e_3+1}=27$.
For $3\nmid a$, put $r=v_3(a^2-1)-1$ and $d=\min\{r,2\}$.  Then
\[
\begin{array}{c|c|c}
a\!\!\!\pmod{27} & r & (k_{3,1},\dots,k_{3,8})\ \text{(nonzero entries)}\\ \hline
\pm1\!\!\pmod3\ \not\equiv \pm1\!\!\pmod9 & 0 & \text{all }0\\
\pm1\!\!\pmod9\ \not\equiv \pm1\!\!\pmod{27} & 1 & k_{3,6}=k_{3,7}=k_{3,8}=1\\
\pm1\!\!\pmod{27} & \ge2 & k_{3,6}=k_{3,7}=1,\ k_{3,8}=2
\end{array}
\]
again matching the explicit basis in \cite{JKS21}.  
\end{example} 

\begin{remark}
Separating the $p\mid n$ denominator profile from the factors $C_m(a)$ supported on primes
dividing $a$ is standard in modern explicit constructions and is exactly the viewpoint in
\cite{JKS21}.  For the restriction to primes $p\mid n$ in denominators see Remark~\ref{fact:denom}, and for $n^2$–periodicity in small degrees see \cite{GR17}.  For the family $x^{p^r}-a$ and variants (with or without $(a,n)=1$), see \cite{Kch24}. 
\end{remark}

By Theorem~\ref{thm:minimal-modulus}, $\mathbf S_p(a)$ depends only on $a\bmod p^{e+1}$; there are only
$p^{e+1}$ residue classes. For each class, the corresponding exponents $k_{p,m}$ and reduced corrections
$\beta_m\bmod p^{k_{p,m}}$ can be read off from the explicit construction of
Jakhar--Khanduja--Sangwan \cite[Theorem~1.6 and Theorem~1.8]{JKS21}. It follows that there are only finitely
many possible $p$--local shapes, and the map
\[
a\ \mapsto\ \mathbf S_p(a)
\]
factors through $a\bmod p^{e+1}$.

As a consequence, for each fixed $n$ there exists a finite table $\mathcal T(n)$ indexed by $a\bmod M(n)$ such that, for every parameter $a$ satisfying the standing hypotheses (in particular Hypothesis~\textup{(\textbf{H})} and the irreducibility of $x^n-a$), the shape $\mathbf S(a)$ of an integral basis of $\mathcal O_{K_a}$ is read off from $\mathcal T(n)$ at the entry $[a\bmod M(n)]$. The table is obtained by concatenating the local classes.

\begin{theorem}[Computation table]\label{thm:precompute}
Fix $n$. There exists a finite lookup table $\mathcal T(n)$, indexed by the residue classes $a\bmod M(n)$, such that for every parameter $a$ satisfying the standing hypotheses the shape of an integral basis of $\mathcal O_{K_a}$ is read off directly from the entry $[a\bmod M(n)]\in\mathcal T(n)$.  The table can be filled purely locally, prime–by–prime at $p\mid n$, using the explicit local recipes in \cite{JKS21} for residue classes modulo $p^{e_p+1}$.

\end{theorem}

\begin{proof}
By Theorem~\ref{thm:minimal-modulus} the $p$–part of the shape depends only on
$a\bmod p^{e_p+1}$.  Thus $\mathcal T(n)$ is the Cartesian product of the local tables at
$p\mid n$ listing, for each residue class modulo $p^{e_p+1}$, the exponents $k_{p,m}$ and the
reduced corrections $\beta_m\bmod p^{k_{p,m}}$; these are exactly what the explicit construction in \cite[Theorem~1.6 and Theorem~1.8]{JKS21} provides. The global shape is obtained by Chinese remaindering of the local choices.
\end{proof}

We can develop relative version for the periodicity as follows. 

\begin{proposition}[Relative periodicity over tame bases]\label{prop:relative}
Let $L/\Bbb Q$ be a number field with $\gcd(n,\Delta_L)=1$, and let $K_{a,L}:=L(\sqrt[n]{a})$
with $a\in\Bbb Z$ $n$th–power–free.  For each rational $p\mid n$ and each prime $\mathfrak p\mid p$
in $\mathcal O_L$, the $\mathfrak p$–local shape of an $\mathcal O_L$–integral basis of
$K_{a,L}/L$ is determined by $a\bmod p^{e_p+1}$.  Consequently the global \emph{relative} shape is periodic modulo $M(n)$.
\end{proposition}

\begin{proof}
Since $\gcd(n,\Delta_L)=1$, every rational prime $p\mid n$ is unramified in $L$, so the local behavior at each prime $\mathfrak p\mid p$ is tame over $L$. The local denominator exponents and reduced correction terms in an $\mathcal O_L$--integral basis of $L(\sqrt[n]{a})/L$ at $\mathfrak p\mid p$ are governed by the same invariant $d_p(a)$ as in the absolute
case (cf.\ \cite{JKS21}), hence are determined by $a\bmod p^{e_p+1}$. Gluing over all $\mathfrak p\mid p$ and all $p\mid n$ yields periodicity modulo $M(n)$.
\end{proof}

The proof above mirrors the Kummer–cyclotomic situation where the decisive congruence lives modulo $(1-\zeta_p)^2$, cf.\ \cite[Theorem.\,1.1 and \S6]{Smi21}.

\begin{example}[Relative periodicity over a quadratic tame base]\label{ex:relative-sextic}
Let $n=6$ and $L=\Bbb Q(\sqrt{13})$. Then $\gcd(n,\Delta_L)=1$ (here $\Delta_L=13$), so $2,3$ are unramified in $L$. In this field, $2$ is inert and $3$ splits:
\[
(2)=\mathfrak p_2,\qquad (3)=\mathfrak p_{3,1}\,\mathfrak p_{3,2}.
\]
Consider $K_{a,L}:=L(\sqrt[6]{a})$ with $a\in\Bbb Z$ $6$th–power–free. By Proposition~\ref{prop:relative}, for each prime $\mathfrak p\mid 6$ the $\mathfrak p$–local shape of an $\mathcal O_L$–integral basis is determined by the class of $a$ modulo $p^{e_p+1}$ with $p\in\{2,3\}$ and $e_2=e_3=1$. Hence the \emph{global relative} shape is periodic modulo
\[
M(6)=2^{\,2}\,3^{\,2}=36,
\]
and this modulus is minimal.

\smallskip
\noindent\emph{Explicit instance.} Take $a\equiv 1\pmod 4$ and $a\equiv 1\pmod 9$ (e.g.\ $a=73$).
\begin{itemize}
\item At $\mathfrak p_2\mid 2$ (inert): $d_2=1$, so
\[
(k_{\mathfrak p_2,1},k_{\mathfrak p_2,2},k_{\mathfrak p_2,3},k_{\mathfrak p_2,4},k_{\mathfrak p_2,5})=(0,0,1,1,1),
\]
with reduced $2$–adic corrections (mod $2$) as in the absolute sextic case:
\[
\beta_{3,\mathfrak p_2}\equiv 1,\quad \beta_{4,\mathfrak p_2}\equiv \theta,\quad
\beta_{5,\mathfrak p_2}\equiv \theta^2\pmod{2}.
\]
\item At each split prime $\mathfrak p_{3,i}\mid 3$ ($i=1,2$): $d_3=1$, so
\[
(k_{\mathfrak p_{3,i},1},k_{\mathfrak p_{3,i},2},k_{\mathfrak p_{3,i},3},k_{\mathfrak p_{3,i},4},k_{\mathfrak p_{3,i},5})=(0,0,0,1,1),
\]
with the same $3$–adic reductions (mod $3$) for both $\mathfrak p_{3,1}$ and $\mathfrak p_{3,2}$:
\[
\beta_{4,\mathfrak p_{3,i}}\equiv 1+\theta^2,\qquad
\beta_{5,\mathfrak p_{3,i}}\equiv \theta+\theta^3\pmod{3}.
\]
\end{itemize}
Thus, over $L$, the only effect of the base field is \emph{how many} primes lie above each $p\mid 6$
(inert vs.\ split); the local exponents and reduced corrections at each $\mathfrak p$ are still decided entirely by
\[
a\bmod 4\quad\text{and}\quad a\bmod 9.
\]
In particular, the relative shape for $K_{a,L}/L$ depends only on $a\bmod 36$, and the four CRT cases from Example~\ref{ex:Sextic-shape} carry over verbatim, with the $3$–adic row duplicated on the two primes $\mathfrak p_{3,1},\mathfrak p_{3,2}$.
\end{example}

\section{Applications of the minimal modulus}\label{sec:applications}

As in \S\ref{sec:defs}, when we write an integral basis in the $\{1,\theta,\dots,\theta^{n-1}\}$–coordinates
we use the Jakhar--Khanduja--Sangwan form
\[
\Bigl\{\,1,\ \frac{\theta^m+\beta_m}{\;C_m(a)\,\prod_{p\mid n}p^{\,k_{p,m}}\;}\ :\ 1\le m\le n-1\Bigr\},
\]
where $C_m(a)$ is supported only at primes dividing $a$, while the $p$–denominator exponents $k_{p,m}\ge0$ and the linear corrections $\beta_m$ are determined locally at each $p\mid n$ from a finite amount of $p$–adic information on $a$.

\subsection{Sharp local index tests at primes \(p\mid n\)}\label{subsec:monogeneity}

We isolate two quick consequences for power–basis tests that fall directly out of the local
classification.

\begin{proposition}[Wieferich threshold when \(e_p=1\)]\label{prop:wieferich}
Suppose $p\parallel n$ and $v_p(a)=0$ (equivalently, $(a,p)=1$). Then
\[
p\mid[\mathcal O_{K_a}:\Bbb Z[\theta]]
\quad\Longleftrightarrow\quad
a^{p-1}\equiv 1 \pmod{p^2}.
\]
Equivalently, $\Bbb Z[\theta]$ is $p$--maximal if and only if $a^{p-1}\not\equiv 1\pmod{p^2}$.
\end{proposition}

\begin{proof}
When $p\parallel n$ (so $e_p=1$) and $v_p(a)=0$, we have
\[
d_p(a)=\min\{\,v_p(a^{p-1}-1)-1,\,1\,\}\in\{0,1\}.
\]
If $a^{p-1}\not\equiv 1\pmod{p^2}$ then $v_p(a^{p-1}-1)=1$, so $d_p(a)=0$ and the explicit basis in
Theorem~\ref{fact:JKS-form} has no $p$--power denominators; equivalently,
$v_p([\mathcal O_{K_a}:\Bbb Z[\theta]])=0$.
If $a^{p-1}\equiv 1\pmod{p^2}$ then $v_p(a^{p-1}-1)\ge 2$, so $d_p(a)=1$ and the explicit basis has a
nontrivial $p$--denominator, hence $v_p([\mathcal O_{K_a}:\Bbb Z[\theta]])>0$.
\end{proof}

In this paper we work in the fixed coordinates $\{1,\theta,\dots,\theta^{n-1}\}$ and in the normalized form where all factors supported on primes dividing $a$ are absorbed into the numerators (via the $C_m(a)$).  In this normalization, the only primes that can occur in the denominators are those dividing $n$ (see Remark~\ref{fact:denom}).  Accordingly, Proposition~\ref{prop:wieferich} concerns the $p$--part of the index $[\mathcal O_{K_a}:\Bbb Z[\theta]]$ at primes $p\mid n$; it is not meant to be a criterion involving primes dividing $a$ (indeed we assume $v_p(a)=0$ here).

\begin{remark}[Index vs.\ monogeneity]\label{rem:index-vs-monogeneity}
Proposition~\ref{prop:wieferich} is a criterion for the $p$--adic contribution to the index
$[\mathcal O_{K_a}:\Bbb Z[\theta]]$, i.e.\ for $p$--maximality of the power order $\Bbb Z[\theta]$. It does \emph{not} decide whether $K_a$ is monogenic (whether there exists some $\alpha\in\mathcal O_{K_a}$ with $\mathcal O_{K_a}=\Bbb Z[\alpha]$).
\end{remark}

\begin{remark}[ABC and the Wieferich threshold at $p\parallel n$]
For $p\parallel n$ and $(a,p)=1$, our order--$1$ local shape at $p$ is decided by the
Wieferich congruence $a^{p-1}\equiv1\pmod{p^2}$ (Proposition~\ref{prop:wieferich}). 
Assuming the $abc$ conjecture, Silverman proved that for any fixed $A\ge2$ there are 
$\gg\log x$ primes $p\le x$ with $A^{p-1}\not\equiv1\pmod{p^2}$; in particular there are 
infinitely many non--Wieferich primes for base $A$ \cite{Silverman88}. Consequently, for a fixed parameter $A$ and varying $n$, the order--$1$ local shape at $p\mid n$ is the $p$--regular (no-denominator) branch for infinitely many $p$. 
Analogous non--Wieferich abundance over number fields under $abc$ is known as well (see, e.g., Ding \cite{Ding19}).
\end{remark}

\begin{proposition}[Local sharpness at every \(p^e\parallel n\)]\label{prop:local-sharpness}
Fix $p^e\parallel n$ with $e\ge1$ and assume $(a,p)=1$.  There are integers
\[
a\equiv 1+p^{e}u\not\equiv 1\pmod{p^{e+1}}\quad\text{and}\quad
a'\equiv 1+p^{e+1}u'\pmod{p^{e+2}}
\]
($u,u'\in\Bbb Z$, $p\nmid uu'$) such that the $p$–parts of the shapes of $K_a$ and $K_{a'}$
differ.  Equivalently, the local modulus $p^{e+1}$ is minimal.
\end{proposition}

\begin{proof}
For odd $p$, the LTE lemma gives
$v_p\bigl((1+p^{e}u)^{p-1}-1\bigr)=e$, hence $r_p(a)=e-1$, while $r_p(a')\ge e$.
Thus $d_p(a)=e-1<e=d_p(a')$ and the local denominator pattern changes with $a\bmod
p^{e+1}$.  For $p=2$ the same conclusion holds with $r_2(a)=v_2(a-1)-1$ by direct
expansion.  The resulting change of $k_{p,m}$ (and of the reduced corrections
$\beta_m\bmod p^{k_{p,m}}$) is built into the Jakhar--Khanduja--Sangwan explicit bases.
\end{proof}

Write \(p^e\parallel n\) and \(n_p:=n/p^e\).  In the unit branch \(p\nmid a\), put
\[
r_p(a):=v_p(a^{p-1}-1)-1,\qquad
t_p(a):=\min\{r_p(a),e\}.
\]
By the Jakhar--Khanduja--Sangwan discriminant formula, specialised at \(p^e\parallel n\),
one has
\begin{equation}\label{eq:disc-exp-t}
v_p(d_{K_a})
=
ne-2n_p\sum_{j=1}^{t_p(a)}p^{e-j},
\end{equation}
with the convention that the empty sum is \(0\).

\begin{proposition}
\label{prop:disc-jump}
Let \(p^e\parallel n\), \(n_p=n/p^e\), and let \(0\le t<e\).  Then there exist
\(n\)th-power-free integers \(a,a'\), both coprime to \(n\), satisfying the standing
hypotheses, such that \( t_p(a)=t,\qquad t_p(a')=t+1,\) and \(a\equiv a'\pmod {p^{t+1}}\), \(a\not\equiv a'\pmod {p^{t+2}}.\) For any such pair,
\[
v_p(d_{K_{a'}})
=
v_p(d_{K_a})-2n_p p^{e-(t+1)}.
\]
In particular, the final jump \(t=e-1\to e\) is realised by two parameters congruent
modulo \(p^e\) but not modulo \(p^{e+1}\).  
\end{proposition}

\begin{proof}
If \(t_p(a)=t\) and \(t_p(a')=t+1\), then \eqref{eq:disc-exp-t} gives
\[
v_p(d_{K_{a'}})-v_p(d_{K_a})
=
-2n_p
\left(
\sum_{j=1}^{t+1}p^{e-j}
-
\sum_{j=1}^{t}p^{e-j}
\right)
=
-2n_p p^{e-(t+1)}.
\]

It remains to produce such \(a,a'\).  Choose residue classes
\[
A\equiv 1+p^{t+1}\pmod {p^{e+2}},
\qquad
A'\equiv 1+p^{t+2}\pmod {p^{e+2}},
\]
and impose, for every prime \(\ell\mid n\), \(\ell\ne p\),
\[
A\equiv A'\equiv 1\pmod \ell .
\]
By the Chinese Remainder Theorem these define primitive classes modulo
\[
Q:=p^{e+2}\prod_{\substack{\ell\mid n\\ \ell\ne p}}\ell .
\]
By Dirichlet's theorem, choose rational primes
\[
a\equiv A\pmod Q,\qquad a'\equiv A'\pmod Q.
\]
Then \(a,a'\) are squarefree, coprime to \(n\), hence \(n\)th-power-free and satisfying
Hypothesis \textup{(\(\mathbf H\))}; moreover \(x^n-a\) and \(x^n-a'\) are irreducible by
Eisenstein. From the construction,
\[
v_p(a-1)=t+1,\qquad v_p(a'-1)=t+2.
\]
For odd \(p\), since \(a\equiv a'\equiv1\pmod p\), LTE gives
\[
v_p(a^{p-1}-1)=t+1,\qquad
v_p((a')^{p-1}-1)=t+2.
\]
For \(p=2\), the same identities are immediate because \(p-1=1\).  Hence
\[
r_p(a)=t,\qquad r_p(a')=t+1,
\]
and since \(0\le t<e\),
\[
t_p(a)=t,\qquad t_p(a')=t+1.
\]
Also
\[
a-a'\equiv p^{t+1}-p^{t+2}=p^{t+1}(1-p)\pmod {p^{t+2}},
\]
so
\[
a\equiv a'\pmod {p^{t+1}},
\qquad
a\not\equiv a'\pmod {p^{t+2}}.
\]

Finally, if \(p\nmid a,a'\) and \(a\equiv a'\pmod {p^{e+1}}\), then
\[
a^{p-1}\equiv (a')^{p-1}\pmod {p^{e+1}},
\]
so the truncated valuations defining \(t_p(a)\) and \(t_p(a')\) agree.  Therefore
\(t_p(a)=t_p(a')\), and \eqref{eq:disc-exp-t} gives
\(v_p(d_{K_a})=v_p(d_{K_{a'}})\).  The case \(t=e-1\) above shows that this constancy
already fails modulo \(p^e\).
\end{proof}

\subsection{Counting and density of shape classes}\label{subsec:counting}
By Theorem~\ref{thm:minimal-modulus}, for fixed $n$ and for parameters satisfying the standing hypotheses, the global shape $\mathbf S(a)$ depends only on the residue class of $a$ modulo $M(n)$. We first isolate the analytic input by counting all $n$th--power--free integers in prescribed arithmetic progressions. This auxiliary count does not impose the irreducibility of $x^n-a$. When we return to degree-$n$ pure fields, we restore Hypothesis~\textup{(\textbf{H})}, the condition $a\neq\pm1$, and irreducibility; Lemma~\ref{lem:reducible-binomial-parameters} shows that the excluded reducible binomial parameters contribute only $O_n(X^{1/2})$ and hence do not change the natural densities.

Recall that a nonzero integer $a$ is \emph{$n$th-power-free} if no prime $p$ satisfies $p^n\mid a$. Equivalently, in the prime factorization $|a|=\prod_p p^{\alpha_p}$ one has $\alpha_p<n$ for every $p$. The sign of $a$ plays no role for $n$th-power-freeness.

We use a standard counting function for $n$th--power--free integers in arithmetic progressions. Fix integers $q\ge1$ and $r\ (\mathrm{mod}\ q)$. We define
\[
\mathcal F_n(X;q,r)
:=\#\bigl\{\,a\in\Bbb Z:\ 0<|a|\le X,\ a\equiv r\ (\mathrm{mod}\ q),\ a\ \text{is $n$th-power-free}\,\bigr\}.
\]
Thus $\mathcal F_n(X;q,r)$ counts $n$th--power--free integers in the arithmetic progression $a\equiv r\pmod q$, up to height $X$ (counting both positive and negative integers; this explains the factor $2X$ that will appear in main terms).  If $\mathcal F_n(X;q,r)$ grows linearly with $X$, then the leading constant is exactly the \emph{natural density} (with respect to height $|a|$) of $n$th--power--free integers in that progression.  Concretely, whenever the limit exists,
\[
\delta_n(q,r)\ :=\ \lim_{X\to\infty}\frac{\mathcal F_n(X;q,r)}{2X}
\]
is the proportion of integers (ordered by $|a|$) that lie in the class $r\bmod q$ and are
$n$th--power--free.  Theorem~\ref{thm:equidistribution} computes this density and makes explicit how it depends on $r$.

Not every residue class modulo $q$ can contain infinitely many $n$th--power--free integers. For example, if $q$ contains a prime power $p^\alpha$ with $\alpha\ge n$ and $r\equiv0\pmod{p^\alpha}$, then every $a\equiv r\pmod q$ is divisible by $p^n$, hence cannot be $n$th--power--free (except possibly $a=0$, which is excluded by $n$th--power--freeness).  The following definitions isolate precisely this obstruction

\begin{definition}[Admissibility notion]\label{def:admissibilities}
Fix $q\ge1$ and $n\ge2$, and let $r$ be a residue class modulo $q$.  If
$p^\alpha\parallel q$, define the $p$--adic valuation of $r$ modulo $q$ by
\[
 v_{p,q}(r):=\max\{\,0\le s\le \alpha:\ r\equiv0\pmod{p^s}\,\}.
\]
Equivalently, if $\widetilde r$ is any integer representative of the class $r$, then
\[
 v_{p,q}(r)=\min\{v_p(\widetilde r),\alpha\},
\]
with the convention $v_p(0)=+\infty$.  Thus $v_{p,q}(r)$ is independent of the chosen
representative.  When $q$ is fixed, we write simply $v_p(r)$ for $v_{p,q}(r)$.
Moreover, if $g\mid q$, then the notation $g\mid r$ means $r\equiv0\pmod g$; this is again
well defined for a residue class modulo $q$.
\begin{itemize}
\item[(a)] We say that $r$ is \emph{$n$--admissible} (weak form) if for every
$p^\alpha\parallel q$ with $\alpha\ge n$ we have $v_p(r)<n$.  Equivalently, the progression
$a\equiv r\pmod q$ does not force the divisibility $p^n\mid a$ for any prime $p$.
\item[(b)] We say that $r$ is \emph{strictly $n$--admissible} if for every $p^\alpha\parallel q$
we have $v_p(r)<\min\{\alpha,n\}$.  This is stronger than (a) and is exactly the condition
that makes the main term below independent of $r$.
\end{itemize}
\end{definition}

\begin{remark}[Terminology]\label{rem:admissibility-terminology}
Definition~\ref{def:admissibilities} concerns $n$th--power--free integers in a general arithmetic progression
modulo $q$. Later, for the specific modulus $q=M(n)$ coming from the local theory, we will often restrict to the
smaller set of residue classes with $v_p(r)\le 1$ for every $p\mid n$; to avoid confusion we will call these
\emph{shape--admissible} classes (Definition~\ref{def:shape-admissible} below).
\end{remark}

Throughout the rest of this subsection, $\mathbf 1_{\mathcal E}$ denotes the indicator of the condition
$\mathcal E$: it is equal to $1$ if $\mathcal E$ holds and to $0$ otherwise.  In particular, since the integer
$g$ occurring below always divides $q$, the symbol $\mathbf 1_{\{g\mid r\}}$ means
\[
\mathbf 1_{\{g\mid r\}}=
\begin{cases}
1,&\text{if }r\equiv0\pmod g,\\
0,&\text{otherwise.}
\end{cases}
\]

For a nonzero integer $a$, the indicator of $n$th--power--freeness is
\[
\mathbf 1_{\text{$n$--free}}(a)
:=
\begin{cases}
1,&\text{if $p^n\nmid a$ for every prime $p$,}\\
0,&\text{otherwise.}
\end{cases}
\]
A classical identity (see, e.g., \cite[Ch.\,III]{Ten15}) is
\begin{equation}\label{eq:mobius-nfree}
\mathbf 1_{\text{$n$--free}}(a)=\sum_{d^n\mid a}\mu(d).
\end{equation}
Here $\mu$ is the M\"obius function.  One quick justification is to factor the right-hand side as an Euler product: since $\mu(d)=0$ unless $d$ is squarefree, we have
\[
\sum_{d^n\mid a}\mu(d)
=\prod_{p}\Bigl(1+\mu(p)\,\mathbf 1_{\{p^n\mid a\}}\Bigr)
=\prod_{p}\Bigl(1-\mathbf 1_{\{p^n\mid a\}}\Bigr),
\]
which equals $1$ if no prime power $p^n$ divides $a$, and equals $0$ otherwise, exactly as desired. Our counting function $\mathcal F_n(X;q,r)$ ignores $a=0$ anyway, since $0$ is not $n$th--power--free.

We now evaluate $\mathcal F_n(X;q,r)$, keeping track of the precise dependence on $r$.
For the estimate below, the notation
\[
E_n(X;q,r)=O_n\!\bigl(X^{1/n}\bigr)
\]
means that there is a constant $C_n>0$, depending only on $n$, such that
$|E_n(X;q,r)|\le C_nX^{1/n}$ for every $X\ge1$, every $q\ge1$, and every residue class
$r\pmod q$. Thus the error estimate is uniform in $q$ and $r$.

\begin{theorem}[Equidistribution with explicit $r$-dependence]\label{thm:equidistribution}
Fix $n\ge2$. For every $q\ge1$, every residue class $r\pmod q$, and every real $X\ge1$,
\begin{equation}\label{eq:Fn-mainterm}
\begin{aligned}
\mathcal F_n(X;q,r)
={}&
\frac{2X}{q}\cdot\frac{1}{\zeta(n)}
\prod_{p\mid q}\Bigl(1-\frac1{p^{\,n}}\Bigr)^{-1}\\
&\times
\prod_{p^\alpha\parallel q}
\Bigl(1-\mathbf 1_{\{v_p(r)\ge \min(\alpha,n)\}}
      p^{\min(\alpha,n)-n}\Bigr)
+E_n(X;q,r),
\end{aligned}
\end{equation}
where
\begin{equation}\label{eq:Fn-error-explicit}
\bigl|E_n(X;q,r)\bigr|
\le
\left(2+\frac{2n}{n-1}\right)X^{1/n}.
\end{equation}
In particular, $E_n(X;q,r)=O_n(X^{1/n})$, uniformly in $q$ and $r$.
For fixed $q$ and $r$, formula \eqref{eq:Fn-mainterm} is therefore an asymptotic formula as
$X\to\infty$, because $X^{1/n}=o(X)$.

Here $v_p(r)=v_{p,q}(r)$ is the residue-class valuation from Definition~\ref{def:admissibilities}.
In particular:
\begin{itemize}
\item If $r$ is \emph{not} $n$--admissible, then $\mathcal F_n(X;q,r)=0$ for every $X$.
\item If $r$ is \emph{strictly $n$--admissible}, then the $r$-dependent product in
\eqref{eq:Fn-mainterm} equals $1$, and
\[
\mathcal F_n(X;q,r)=
\frac{2X}{q}\cdot\frac{1}{\zeta(n)}
\prod_{p\mid q}\Bigl(1-\frac1{p^{\,n}}\Bigr)^{-1}
+O_n\!\bigl(X^{1/n}\bigr),
\]
independently of $r$.
\end{itemize}
\end{theorem}

\begin{proof}
Put $Y=X^{1/n}$. By definition and by \eqref{eq:mobius-nfree},
\[
\mathcal F_n(X;q,r)
=\sum_{\substack{0<|a|\le X\\ a\equiv r\ (\mathrm{mod}\ q)}}
 \sum_{d^n\mid a}\mu(d).
\]
If $d^n\mid a$ and $0<|a|\le X$, then $d\le Y$. Hence the double sum is finite, and we may
interchange the order of summation:
\begin{equation}\label{eq:Fn-swap}
\mathcal F_n(X;q,r)
=\sum_{1\le d\le Y}\mu(d)\,
\#\bigl\{\,0<|a|\le X:\ a\equiv r\pmod q,\ d^n\mid a\,\bigr\}.
\end{equation}

Fix such a $d$ and put
\[
g_d:=\gcd(q,d^n),\qquad
L_d:=\operatorname{lcm}(q,d^n)=\frac{qd^n}{g_d}.
\]
The simultaneous congruences
\[
a\equiv r\pmod q,
\qquad a\equiv0\pmod{d^n}
\]
are compatible if and only if $r\equiv0\pmod{g_d}$, i.e. if and only if $g_d\mid r$ in the
residue-class sense of Definition~\ref{def:admissibilities}. If they are compatible, their
solutions form one residue class modulo $L_d$.

We now make the counting error completely explicit. For any modulus $L\ge1$ and any residue
class $c\pmod L$, choose an integer representative $c_0$. The relevant count is
\[
\left\lfloor\frac{X-c_0}{L}\right\rfloor
-\left\lceil\frac{-X-c_0}{L}\right\rceil+1,
\]
and therefore differs from $2X/L$ by at most $1$. Removing the integer $0$ changes this count
by at most one. Consequently, for every $d$ there is a real number $\varepsilon_d(X;q,r)$ with
$|\varepsilon_d(X;q,r)|\le2$ such that
\begin{equation}\label{eq:single-progression-count}
\begin{aligned}
&\#\bigl\{\,0<|a|\le X:\ a\equiv r\pmod q,\ d^n\mid a\,\bigr\}\\
&\qquad=
\mathbf 1_{\{g_d\mid r\}}
\left(\frac{2X}{L_d}+\varepsilon_d(X;q,r)\right)\\
&\qquad=
\mathbf 1_{\{g_d\mid r\}}
\left(\frac{2X}{q}\frac{g_d}{d^n}+\varepsilon_d(X;q,r)\right).
\end{aligned}
\end{equation}
The second equality in \eqref{eq:single-progression-count} is exact, since
$L_d=qd^n/g_d$; no term is absorbed into the error.

Substituting \eqref{eq:single-progression-count} into \eqref{eq:Fn-swap} gives
\begin{equation}\label{eq:Fn-SX}
\mathcal F_n(X;q,r)=\frac{2X}{q}S_X(q,r)+E_1(X;q,r),
\end{equation}
where
\[
S_X(q,r):=
\sum_{1\le d\le Y}
\mu(d)\frac{\gcd(d^n,q)}{d^n}\,
\mathbf 1_{\{\gcd(d^n,q)\mid r\}}
\]
and
\[
E_1(X;q,r):=
\sum_{1\le d\le Y}
\mu(d)\mathbf 1_{\{g_d\mid r\}}\varepsilon_d(X;q,r).
\]
Since $|\mu(d)|\le1$ and $|\varepsilon_d(X;q,r)|\le2$,
\begin{equation}\label{eq:E1-bound}
|E_1(X;q,r)|\le2\lfloor Y\rfloor\le2Y.
\end{equation}
The constant $2$ in this estimate is absolute; in particular, it is independent of $q$ and $r$.

Consider now the absolutely convergent series
\[
S(q,r):=
\sum_{d\ge1}
\mu(d)\frac{\gcd(d^n,q)}{d^n}\,
\mathbf 1_{\{\gcd(d^n,q)\mid r\}}.
\]
Indeed, $\gcd(d^n,q)\le q$ and $n\ge2$. Moreover,
\[
|S(q,r)-S_X(q,r)|
\le q\sum_{d>Y}\frac1{d^n}.
\]
Since $Y\ge1$ and $t\mapsto t^{-n}$ is decreasing,
\[
\sum_{d>Y}\frac1{d^n}
\le Y^{-n}+\int_Y^\infty t^{-n}\,dt
\le \frac{n}{n-1}Y^{1-n}.
\]
Therefore
\begin{equation}\label{eq:tail-bound}
\left|\frac{2X}{q}\bigl(S(q,r)-S_X(q,r)\bigr)\right|
\le\frac{2n}{n-1}XY^{1-n}
=\frac{2n}{n-1}X^{1/n}.
\end{equation}
Combining \eqref{eq:Fn-SX}, \eqref{eq:E1-bound}, and \eqref{eq:tail-bound}, we obtain
\begin{equation}\label{eq:Fn-S}
\mathcal F_n(X;q,r)=\frac{2X}{q}S(q,r)+E_n(X;q,r),
\end{equation}
with the explicit bound \eqref{eq:Fn-error-explicit}. This also explains why the error depends
only on $n$: the factor $q$ in the tail estimate cancels against the prefactor $1/q$.

It remains to compute $S(q,r)$. The summand is multiplicative in $d$. Indeed, if
$\gcd(d_1,d_2)=1$, then
\[
\gcd((d_1d_2)^n,q)=\gcd(d_1^n,q)\gcd(d_2^n,q),
\]
the two factors on the right are coprime, and hence
\[
\mathbf 1_{\{\gcd((d_1d_2)^n,q)\mid r\}}
=
\mathbf 1_{\{\gcd(d_1^n,q)\mid r\}}
\mathbf 1_{\{\gcd(d_2^n,q)\mid r\}}.
\]
Absolute convergence therefore permits an Euler-product factorisation. Since $\mu(d)$ is
supported on squarefree integers, each prime contributes only the choices $p\nmid d$ and
$p\mid d$.

If $p^\alpha\parallel q$, the local term corresponding to $p\mid d$ is
\[
\mu(p)\,p^{\min(n,\alpha)-n}\,
\mathbf 1_{\{v_p(r)\ge\min(n,\alpha)\}},
\]
so the local factor is
\[
1-
\mathbf 1_{\{v_p(r)\ge\min(n,\alpha)\}}
 p^{\min(n,\alpha)-n}.
\]
If $p\nmid q$, the compatibility condition is automatic and the local factor is
$1+\mu(p)p^{-n}=1-p^{-n}$. Hence
\[
S(q,r)=
\prod_{p\nmid q}(1-p^{-n})
\prod_{p^\alpha\parallel q}
\Bigl(1-
\mathbf 1_{\{v_p(r)\ge\min(\alpha,n)\}}
 p^{\min(\alpha,n)-n}\Bigr).
\]
Finally,
\[
\prod_{p\nmid q}(1-p^{-n})
=
\prod_p(1-p^{-n})\prod_{p\mid q}(1-p^{-n})^{-1}
=
\frac1{\zeta(n)}\prod_{p\mid q}(1-p^{-n})^{-1}.
\]
Substitution into \eqref{eq:Fn-S} proves \eqref{eq:Fn-mainterm}.

If $r$ is not $n$--admissible, then for some $p^\alpha\parallel q$ with $\alpha\ge n$ one has
$v_p(r)\ge n$. Thus every integer $a\equiv r\pmod q$ is divisible by $p^n$, and therefore
$\mathcal F_n(X;q,r)=0$. If $r$ is strictly $n$--admissible, then
$v_p(r)<\min(\alpha,n)$ for every $p^\alpha\parallel q$, so all indicators in the second
product vanish. This proves the stated $r$-independent main term.
\end{proof}

Dividing \eqref{eq:Fn-mainterm} by $2X$ and letting $X\to\infty$ shows that the density $\delta_n(q,r)$ exists for every residue class $r$ and is given by the constant multiplying $\frac{2X}{q}$. This is the auxiliary arithmetic density needed below. For the application to shapes of degree-$n$ pure fields, we shall restrict to shape--admissible classes (which automatically enforce Hypothesis~\textup{(\textbf{H})}) and then remove the sparse parameters for which $x^n-a$ is reducible.

\begin{corollary}\label{cor:nfull}
If $q$ is \emph{$n$–full}, i.e.\ $v_p(q)\ge n$ for every $p\mid q$, then all $n$–admissible residues $r\ (\mathrm{mod}\ q)$ are strictly $n$–admissible, and therefore the $n$th–power–free integers are equidistributed among the $n$–admissible classes with common density
\[
\frac{1}{q}\,\zeta(n)^{-1}\,\prod_{p\mid q}\Bigl(1-\frac1{p^{\,n}}\Bigr)^{-1}.
\]
\end{corollary}

In our applications we take $q=M(n)=\prod_{p^e\parallel n}p^{e+1}$. Since typically $e+1<n$, Corollary~\ref{cor:nfull} need not apply, and the uniform density across \emph{all} admissible classes mod $M(n)$ may fail. Theorem~\ref{thm:equidistribution}, however, still yields the correct main term with the explicit $r$–dependence. For many counting questions it suffices to restrict attention to strictly $n$–admissible classes (e.g.\ those with $v_p(r)\in\{0,1\}$ at $p\mid n$), in which case uniformity is restored.

\begin{definition}[Shape--admissible residue classes modulo $M(n)$]\label{def:shape-admissible}
A residue class $r\ (\mathrm{mod}\ M(n))$ is called \emph{shape--admissible} if
\[
v_p(r)\le 1\qquad\text{for every prime }p\mid n.
\]
\end{definition}

The arithmetic estimate of Theorem~\ref{thm:equidistribution} applies to all $n$th--power--free parameters. To obtain a statement about degree-$n$ pure fields, we must also impose the standing irreducibility condition. The next elementary lemma shows that this changes only the error term, not the density.

\begin{lemma}[Sparse reducible binomial parameters]\label{lem:reducible-binomial-parameters}
Fix $n\ge2$, and put
\[
\mathcal E_n(X):=
\#\bigl\{\,a\in\Bbb Z:\ 0<|a|\le X,\ x^n-a\text{ is reducible over }\Bbb Q\,\bigr\}.
\]
Then
\[
\mathcal E_n(X)\ll_n X^{1/2}\qquad(X\ge1).
\]
\end{lemma}

\begin{proof}
By Capelli's criterion for integer binomials \cite{KoleyReddy20}, if $x^n-a$ is reducible over $\Bbb Q$, then either
\[
a=b^t\quad\text{for some }b\in\Bbb Z\text{ and some divisor }t>1\text{ of }n,
\]
or $4\mid n$ and $a=-4b^4$ for some $b\in\Bbb Z$. In the first case, choosing a prime $\ell\mid t$ gives $\ell\mid n$ and $a=c^\ell$ for some $c\in\Bbb Z$. Hence
\[
\mathcal E_n(X)
\le
\sum_{\substack{\ell\mid n\\ \ell\ \text{prime}}}
\bigl(2\lfloor X^{1/\ell}\rfloor+1\bigr)
+
\mathbf 1_{\{4\mid n\}}\bigl(2\lfloor (X/4)^{1/4}\rfloor+1\bigr)
\ll_n X^{1/2}.
\]
\end{proof}

Let $\mathcal A_n$ denote the set of integers $a\neq0,\pm1$ that are $n$th--power--free, for which $(n,a)$ satisfies Hypothesis~\textup{(\textbf{H})}, and for which $x^n-a$ is irreducible over $\Bbb Q$. For $a\in\mathcal A_n$, Theorem~\ref{thm:minimal-modulus} shows that the shape map
\[
\mathsf{Sh}_n:\ \mathcal A_n\longrightarrow\mathcal S_n,
\qquad a\longmapsto\mathbf S(a),
\]
factors through $a\bmod M(n)$. Fix any union
$\mathcal R\subseteq\Bbb Z/M(n)\Bbb Z$ of shape--admissible classes and define
\[
\mathcal N_\mathcal R(X)
:=
\#\bigl\{\,a\in\mathcal A_n:\ |a|\le X,\ a\bmod M(n)\in\mathcal R\,\bigr\}.
\]
This is a count of defining parameters $a$, not of isomorphism classes of number fields.

\begin{theorem}[Equidistribution and density of pure-field parameters]\label{thm:shape-density}
For every fixed $n\ge3$ and every union $\mathcal R$ of shape--admissible classes modulo $M(n)$,
\[
\mathcal N_\mathcal R(X)
=
\frac{2X}{M(n)}\,
\frac{1}{\zeta(n)}
\prod_{p\mid n}\Bigl(1-\frac1{p^{\,n}}\Bigr)^{-1}
\#\mathcal R
+O_n(X^{1/2}).
\]
The implied constant is uniform in $\mathcal R$. In particular, the natural density of admissible degree-$n$ pure-field parameters $a$ with $a\bmod M(n)\in\mathcal R$ exists and equals
\[
\delta(\mathcal R)
=
\frac{\#\mathcal R}{M(n)}\cdot
\frac{1}{\zeta(n)}
\prod_{p\mid n}\Bigl(1-\frac1{p^{\,n}}\Bigr)^{-1}.
\]
\end{theorem}

\begin{proof}
Let $q=M(n)=\prod_{p^e\parallel n}p^{e+1}$. If $r$ is shape--admissible, then for every $p^e\parallel n$ one has $v_p(r)\le1$. Since $1<e+1$, every integer $a\equiv r\pmod{M(n)}$ satisfies
\[
v_p(a)=v_p(r)\in\{0,1\}.
\]
Thus every $n$th--power--free integer in such a class automatically satisfies Hypothesis~\textup{(\textbf{H})}. Moreover,
\[
v_p(r)\le1<\min\{e+1,n\},
\]
so every shape--admissible class is strictly $n$--admissible in the sense of Definition~\ref{def:admissibilities}.

Let
\[
\mathcal F_\mathcal R(X):=
\sum_{r\in\mathcal R}\mathcal F_n(X;M(n),r).
\]
Applying Theorem~\ref{thm:equidistribution} to each $r\in\mathcal R$ gives
\[
\mathcal F_\mathcal R(X)
=
\frac{2X}{M(n)}\cdot\frac{1}{\zeta(n)}
\prod_{p\mid n}\Bigl(1-\frac1{p^{\,n}}\Bigr)^{-1}
\#\mathcal R
+O_n\!\bigl(\#\mathcal R\,X^{1/n}\bigr).
\]
Because $\#\mathcal R\le M(n)$ and $n$ is fixed, this error is $O_n(X^{1/n})$. The count $\mathcal N_\mathcal R(X)$ is obtained from $\mathcal F_\mathcal R(X)$ by removing the parameters for which $x^n-a$ is reducible, together with the at most two exceptional values $a=\pm1$. Lemma~\ref{lem:reducible-binomial-parameters} therefore gives
\[
\bigl|\mathcal F_\mathcal R(X)-\mathcal N_\mathcal R(X)\bigr|
\ll_n X^{1/2}.
\]
Since $n\ge3$, one has $X^{1/n}\le X^{1/2}$ for $X\ge1$, and the asserted asymptotic follows. Division by $2X$ and passage to the limit gives the density formula.
\end{proof}

This is the standard Möbius-inversion count for $k$-free numbers in arithmetic progressions (see \cite{Ten15}), specialized to $k=n$ and the modulus $M(n)$ singled out by the local theory. The additional use of Capelli's criterion removes the zero-density set of reducible binomial parameters. The main term factors as a product of local proportions, matching the CRT-gluing of local shapes. For quantitative refinements concerning $k$-free integers in arithmetic progressions, see Nunes \cite{Nun22} and Parry \cite{Par21}.

\begin{remark}[Local factorisation of the count]\label{rem:local-count}
If the chosen set of classes is a product of independent local conditions,
\[
\mathcal R=\prod_{p\mid n}\mathcal R_p
\]
under the Chinese Remainder identification
\[
\mathbb Z/M(n)\mathbb Z\simeq
\prod_{p^e\parallel n}\mathbb Z/p^{e+1}\mathbb Z,
\]
then \(\#\mathcal R=\prod_{p\mid n}\#\mathcal R_p\), and the corresponding
density factorises into the product of the local proportions.
\end{remark}

\begin{example}[Counting and density for $n=4$]\label{ex:count-n4}
Here $M(4)=8$ and the shape--admissible residue classes (in the sense after Theorem~\ref{thm:shape-density}) are those $r\ (\mathrm{mod}\ 8)$ with $v_2(r)\le 1$, i.e.\ $r\in\{1,3,5,7,2,6\}$ (six classes). The three quartic shapes are:
\[
\begin{array}{c|c}
\text{shape }(k_{2,1},k_{2,2},k_{2,3}) & \text{classes }r\ (\mathrm{mod}\ 8)\\ \hline
(0,1,2) & \{1\}\\
(0,1,1) & \{5\}\\
(0,0,0) & \{3,7,2,6\}
\end{array}
\]
By Theorem~\ref{thm:shape-density}, for any union $\mathcal R$ of shape--admissible classes one has
\[
\delta(\mathcal R)=\frac{\#\mathcal R}{8}\cdot \frac{1}{\zeta(4)}\Bigl(1-\frac1{2^4}\Bigr)^{-1}.
\]
Hence the individual parameter densities of these shapes are
\[
\delta_{(0,1,2)}=\delta_{(0,1,1)}=\frac{1}{8}\cdot \frac{1}{\zeta(4)}\Bigl(1-\frac1{2^4}\Bigr)^{-1}
=\frac{12}{\pi^4},\qquad
\delta_{(0,0,0)}=\frac{4}{8}\cdot \frac{1}{\zeta(4)}\Bigl(1-\frac1{2^4}\Bigr)^{-1}
=\frac{48}{\pi^4}.
\]
\end{example}

\begin{example}[Counting and density for $n=6$]\label{ex:count-n6}
Here $M(6)=36$, and shape--admissible classes are those with $v_2(r)\le 1$ and $v_3(r)\le 1$, i.e.\ $r\ (\mathrm{mod}\ 4)\in\{1,2,3\}$ and $r\ (\mathrm{mod}\ 9)\in\{1,2,3,4,5,6,7,8\}$, giving $3\times 8=24$ shape--admissible classes in total.

Write
\[
C_6\ :=\ \frac{1}{\zeta(6)}\Bigl(1-\frac1{2^6}\Bigr)^{-1}\Bigl(1-\frac1{3^6}\Bigr)^{-1}.
\]
Using the CRT rows in Example~\ref{ex:Sextic-shape}, the four sextic shape‑types occur on the following numbers of classes:
\[
\renewcommand{\arraystretch}{1.15}
\begin{array}{l|c|c}
\text{type} & \text{defining congruences} & \#\text{classes mod }36\\
\hline
\text{both }(2,3)\text{ contributions} & a\equiv 1\ (\mathrm{mod}\ 4),\ a\equiv \pm1\ (\mathrm{mod}\ 9) & 2\\
\text{only $2$} & a\equiv 1\ (\mathrm{mod}\ 4),\ a\not\equiv \pm1\ (\mathrm{mod}\ 9) & 6\\
\text{only $3$} & a\not\equiv 1\ (\mathrm{mod}\ 4),\ a\equiv \pm1\ (\mathrm{mod}\ 9) & 4\\
\text{none} & a\not\equiv 1\ (\mathrm{mod}\ 4),\ a\not\equiv \pm1\ (\mathrm{mod}\ 9) & 12
\end{array}
\]
Therefore, by Theorem~\ref{thm:shape-density}, the corresponding parameter densities are
\[
\begin{aligned}
\delta_{\text{both}}
&=\frac{2}{36}\,C_6=\frac{1}{18}\,C_6,
&\qquad
\delta_{\text{only }2}
&=\frac{6}{36}\,C_6=\frac{1}{6}\,C_6,\\
\delta_{\text{only }3}
&=\frac{4}{36}\,C_6=\frac{1}{9}\,C_6,
&\qquad
\delta_{\text{none}}
&=\frac{12}{36}\,C_6=\frac{1}{3}\,C_6.
\end{aligned}
\]
(These four densities sum to $\tfrac{24}{36}C_6$, reflecting that we have restricted to the $24$ shape--admissible classes.)
\end{example}
\bibliographystyle{alpha}
\bibliography{References_final}

\end{document}